%% file: root.tex
\documentclass[12pt, anon]{l4dc2025}

\input{l4dc_2024/packages}

\definecolor{cello}{HTML}{ffe6cc}

\definecolor{mygray}{gray}{0.5}

\newcommand{\ssp}{\mathbb{X}}
\newcommand{\compacts}{\mathbb{K}\ssp}
\newcommand{\ti}{t}
\newcommand{\tf}{t_f}
\newcommand{\intt}{\tau}
\newcommand{\tint}{\mathbb{T}}

\newcommand{\sti}{x}
\newcommand{\tj}{\mathsf{x}}
\newcommand{\sys}{f}

\newcommand{\tg}{\mathcal{T}}
\newcommand{\tgfn}{J_{\tg}}
\newcommand{\ctrlset}{\mathcal{U}}
\newcommand{\distset}{\mathcal{D}}
\newcommand{\csig}{\mathsf{u}}
\newcommand{\dsig}{\mathsf{d}}
\newcommand{\ctrlsetsig}{\mathbb{U} (\ti)}
\newcommand{\distsetsig}{\mathbb{D} (\ti)}
\newcommand{\strd}{\mathfrak{d}}
\newcommand{\stru}{\mathfrak{u}}
\newcommand{\stratset}{\mathfrak{D}(\ti)}
\newcommand{\stratsetu}{\mathfrak{U}(\ti)}
\newcommand{\hjrset}{\mathcal{R}}

\newcommand{\fset}{\mathcal{S}}

\newcommand{\hopf}{\mathcal{H}}

\newcommand{\ham}{H}

\newcommand{\minou}{\min_{u \in \ctrlset}}
\newcommand{\maxod}{\max_{d \in \distset}}

\newcommand{\supod}{\sup_{\strd \in \stratset}}
\newcommand{\infou}{\inf_{\csig \in \ctrlsetsig}}
\newcommand{\supou}{\sup_{\stru \in \stratsetu}}
\newcommand{\infod}{\inf_{\dsig \in \distsetsig}}

\newcommand{\minot}{\min_{\tau \in \tint}}

\newcommand{\sasp}{\tilde{\ssp}}

\newcommand{\maxerr}{\delta^*}
\newcommand{\err}{\varepsilon}
\newcommand{\esig}{\mathbf{\bm{\err}}}
\newcommand{\stre}{\mathfrak{e}}
\newcommand{\stratsete}{\mathfrak{E}(\ti)}
\newcommand{\infoe}{\inf_{\mathfrak{e} \in \mathfrak{E}(\ti)}}
\newcommand{\supoe}{\sup_{\mathfrak{e} \in \mathfrak{E}(\ti)}}

\newcommand{\lsys}{\ell}
\newcommand{\lsyse}{\lsys + \err}
\newcommand{\lintj}{\tj_\lsys}
\newcommand{\eset}{\mathcal{E}}
\newcommand{\esigset}{\mathbb{E} (\ti)}

\newcommand{\hamlerr}{\ham_{\err}}
\newcommand{\loss}{\mathcal{L}}
\newcommand{\losspde}{\loss_{PDE}}
\newcommand{\losslssd}{\loss_{LSS-D}}
\newcommand{\losslssnl}{\loss_{LSS-NL}}

\newcommand{\lossterm}{L}
\newcommand{\losstermpde}{\lossterm_{PDE}}
\newcommand{\losstermls}{\lossterm_{LS}}

\newcommand{\augp}{\lambda}
\newcommand{\augtj}{\tilde{\mathsf{x}}}
\newcommand{\augsys}{\tilde{\sys}_{\augp}}

\newtheorem{Theorem}{Theorem}
\newtheorem{Lemma}{Lemma}

\newtheorem{Definition}{Definition}
\newtheorem{Corollary}{Corollary}

\newtheorem{Remark}{Remark}

\title{
Linear Supervision for Nonlinear, High-Dimensional Neural Control and Differential Games
}

\author{
 \Name{William Sharpless} \Email{wsharpless@ucsd.edu}\\
 \addr Department of Mechanical and Aerospace Engineering, University of California San Diego
 \AND
 \Name{Zeyuan Feng} \Email{zeyuanf@stanford.edu}\\
  \addr Department of Aeronautics and Astronautics, Stanford University
 \AND
 \Name{Somil Bansal} \Email{somil@stanford.edu}\\
 \addr Department of Aeronautics and Astronautics, Stanford University
 \AND
 \Name{Sylvia Herbert} \Email{sherbert@ucsd.edu}\\
 \addr Department of Mechanical and Aerospace Engineering, University of California San Diego%
}

\begin{document}

\maketitle

\thanks{This work is supported by ONR YIP N00014-22-1-2292. The content is solely the responsibility of the authors.}

\begin{abstract}

As the dimension of a system increases, traditional methods for control and differential games rapidly become intractable, making the design of safe autonomous agents challenging in complex or team settings. Deep-learning approaches avoid discretization and yield numerous successes in robotics and autonomy, but at a higher dimensional limit, accuracy falls as sampling becomes less efficient. We propose using rapidly generated \textit{linear} solutions to the partial differential equation (PDE) arising in the problem to accelerate and improve learned value functions for guidance in high-dimensional, \textit{nonlinear} problems. We define two programs that combine supervision of the linear solution with a standard PDE loss. We demonstrate that these programs offer improvements in either speed or accuracy in both a 50-D differential game problem and a 10-D quadrotor control problem.

\end{abstract}

\section{Introduction}

The goal of this work is to solve differential games and optimal controllers for high-dimensional, nonlinear problems. These problems arise in many fields, including multi-agent robotics, medicine, and finance (\cite{Bansal2017}). High-dimensional nonlinearity renders most approaches to autonomy in these domains infeasible, making this a difficult challenge. 
Ultimately, we aim to generate a value function that an agent or team may follow to safely guide itself amidst disturbances or antagonists.

A popular method for the analysis of differential games is to solve its corresponding Hamilton-Jacobi (HJ) partial differential equation (PDE) (\cite{Evans84}). The sub-zero level set of the solution represents the \textit{Backwards Reachable Tube} (BRT) of states that may reach or avoid a target despite any bounded disturbance (\cite{Mitchell05, altarovici2013general}). Moreover, the gradient of the solution defines the corresponding control law. This approach is widely successful (\cite{Huang2011CaptureTheFlag, Chen2015Platooning, jeong2024robots}). However, numerical integration of the PDE with dynamic programming requires a grid, leading to an exponential growth in computation with dimension, becoming infeasible beyond seven.
To overcome this ``curse of dimensionality'', many directions have been explored. Decomposition of the system offers a conservative value when feasible (\cite{chen2018decomposition, he2023efficient}). Multiple methods directly approximate the BRT, largely for linear and polynomial dynamics (\cite{althoff2008reachability, Majumdar2017, Kousik2018RTD, yang2022efficient, kochdumper2023constrained}). These formal methods offer strong guarantees but can suffer from conservativeness.

For convex games with linear time-varying dynamics, several have recognized that a solution to the HJ-PDE is given by the Hopf formula (\cite{darbon2016algorithms, chow2017algorithm, kurzhanski2014dynamics}). By computing the value independently in space and time, this approach offers orders of magnitude in acceleration. Some have used this to approximate the value in nonlinear settings (\cite{Kirchner_2018, lee19iterativehopf, sharpless2023koopman}), but few offer guarantees (\cite{kurzhanski2014dynamics, liu2024hamilton, sharpless2024state}), typically at the cost of being overly conservative.

With the advent of deep learning, many have sought to learn the value function. Typically, a residual based PDE loss (\cite{raissi2017physics, han2018solving}) is used to train the approximation to satisfy the HJ equation (\cite{djeridane2006neural, darbon2020overcoming}). This has similarly been explored in reinforcement learning (\cite{fisac2019bridging, ganai2024hamilton}). One may also learn the dynamic programming procedure itself (\cite{esteve2024finite}). Amongst these, we build on \cite{bansal2020deepreach}, which demonstrated that a sinusoidal architecture (\cite{sitzmann2020implicit}) was well-suited for the problem. Learning the HJ-PDE solution has proved to increase the feasible dimension, however, many of the works demonstrate a nebulous limit beyond which learned solutions deteriorate.

We seek to combine the scalability of linear HJ-PDE's with the efficacy of learned methods to obtain accurate value functions for nonlinear, high-dimensional games. 
% Namely, we propose using linear solutions for semi-supervision to learn nonlinear game value. 
Namely, we introduce supervision with a linear solution to the PDE loss program to learn nonlinear solutions. 
This shifts the learning process from generation to refinement, and results in a faster and more accurate learning process with tighter probabilistic guarantees. 
In summary, \vspace{-0.5em}
%We provide two such methods, a) an augmented game and b) a linear decay, and offer rigorous arguments for each. %We demonstrate these on a 50-D decomposable differential game and a 10-D quadrotor optimal control problem, concluding that our approaches are more accurate and faster, and yield tighter probabilistic guarantees. Explicitly, we

% \begin{figure}[t]
%     \centering
%     \includegraphics[width=0.5\linewidth, height=0.5\linewidth]{fig_sketches/HL_figsk.jpg}    
%     \caption{\textbf{Sketch Figure?} On the left, a depcition of the linear, convex solution which we can learn rapidly by linear supervision with Hopf. On the right, a depiction of the true nonconvex solution. Arrow between denoting linear to nonlinear system. (CUT FOR SPACE?)}
%     \label{fig:TVonly}\vspace{-1em}
% \end{figure}

\begin{enumerate}
    \item We propose two semi-supervision programs for learning high-dimensional, nonlinear value functions. \vspace{-0.5em}
    \item We demonstrate their potential with a 50-D game and 10-D quadrotor control.
    % Here, our methods are more accurate and faster, and yield tighter probabilistic guarantees. \vspace{-1em}
    \item  We provide theoretical insights behind the effectiveness of the proposed programs. \vspace{-0.5em}
\end{enumerate}
% The paper is structured as follows ... 
% Sec. \ref{sec:HJR} formally introduces HJR, BRTs and the DP solution. Sec. \ref{subsec:hopf} introduces the Hopf formula. Sec. \ref{subsec:AntagErr} gives the major results of the work for solving safe envelopes and Sec.\ref{subsec:lim} describes the innate limitations. Sec.~\ref{sec:improvements} offers methods for improving the guarantees, Sec.\ref{sec:det} describes an algorithm for solving the problem.
% % and Sec.~\ref{subsec:fxdctrl} discusses how this work might extend toward fixed-policy validation and improvement
% Finally, Sec.\ref{sec:Demos} demonstrates the theory in a problem of validating fixed policies in a multi-agent pursuit-evasion game using Dubin's car models. 

\section{Problem Statement}\label{sec:prelims}

% Let the evolution of a state in $\ssp \triangleq \mathbb{R}^n$ be defined by a nonlinear function $f$ s.t.
% \begin{eqnarray}
% \dot{\tj} \triangleq \sys(\tj, \csig, \dsig, \intt)
% \label{Dynamics}
% \end{eqnarray}
% where $\tau \in \tint \triangleq[\ti, \tf]$ is the time and $\tj:\tint \to \ssp$ is the unique state trajectory beginning at $\tj(\ti) = x$ arising from the measurable functions $\csig:\tint \to \ctrlset$ and $\dsig:\tint \to \distset$ that represent control and disturbance signals, which are 
% limited to compact and convex sets $\ctrlset \subset \mathbb{R}^{n_u}$ and $\distset \subset \mathbb{R}^{n_d}$. Let \eqref{Dynamics} be Lipschitz continuous in $(\tj, \csig, \dsig)$ and continuous in $\tau$.

% The goal of this work is to guide an autonomous agent in the face of disturbance or adversary. 
We aim to guide the agent state $x \in \ssp \triangleq \mathbb{R}^n$ that evolves by nonlinear dynamics $\sys$ such that we optimize the cost $\tgfn:\ssp \to \mathbb{R}$ of the trajectory at some time $\tau \in \tint \triangleq[\ti, \tf]$. Let the cost be defined such that its sub-zero level set represents a compact, convex target $\tg \subset \ssp$.
While the agent controls the evolution, assume it plays against an opponent who is seeking the opposite score of $\tgfn$ in a zero-sum fashion. 
Hence, if the agent seeks to reach/avoid the target, the opponent seeks to avoid/reach it. 
Note, if the disturbance action is null, the game is an optimal control problem. In this section, let the goal be of reaching the target.

Formally, we seek to solve the differential game,

\begin{equation}
\left\{
\begin{array}{cl}
\underset{\strd \in \stratset}{\operatorname{maximize}} \:\: \underset{\csig \in \ctrlsetsig}{\operatorname{minimize}}  \:\: \underset{\tau \in \tint}{\operatorname{minimize}} & \tgfn (\tj(\tau)) \\

 & \dot{\tj}(\tau) = \sys(\tj(\tau), \csig(\tau), \dsig(\tau), \intt) \\
\text { subject to } & \csig(\tau) \in \ctrlset, \:\: \dsig(\tau) = \strd[\csig](\tau) \in \distset,\\
& \tj(\ti) = x.
\end{array}\right.
\label{prob}
\end{equation}
Here, $\tj:\tint \to \ssp$ is the unique state trajectory beginning at $\tj(\ti) = x$. The trajectory arises from the effect of the agent's control signal $\csig:\tint \to \ctrlset$ and the opponent's disturbance signal $\dsig:\tint \to \distset$, which are in sets of measurable functions $\ctrlsetsig$ and $\distsetsig$ that output actions in the compact, convex sets $\ctrlset \subset \mathbb{R}^{n_u}$ and $\distset \subset \mathbb{R}^{n_d}$ respectively. This disturbance signal is the output of the opponent's strategy $\strd:\ctrlsetsig \to \distsetsig$, which is in the set of non-anticipative strategies $\stratset$ s.t. for $\strd \in \stratset$, $\forall s \in [t,\tau], u(s) = \hat u(s) \implies \strd[u](s) = \strd[\hat u](s)$ \cite{Mitchell05}. Let $\sys$ be Lipschitz continuous in $(\tj, \csig, \dsig)$ and continuous in $\tau$. 

For an initial condition $(\sti,\ti)$, the optimal value of this game is given by
\begin{eqnarray}
\begin{aligned}
V(\sti, \ti) &= \supod \infou \minot \tgfn (\tj(\tau)). \\
\end{aligned}
\label{GameValue}
\end{eqnarray}
Notably, this value describes where the agent may win or lose (\cite{Mitchell05}) as 
\begin{eqnarray}
\begin{aligned}
V(x,t) < 0 \iff x \in \hjrset_{\tg}(\ti) \triangleq \{\sti \in \ssp &\mid \forall \strd \in \stratset, \:\: \exists \csig\in \ctrlsetsig \text{ s.t. } \tj(\tau) \in \tg, \: \tau \in \tint \}.
\end{aligned}
\label{BRT}
\end{eqnarray}
$\hjrset_\tg$ is the BRT that contains all states from which the agent can reach the target despite any disturbance strategy. Moreover, this value also offers the optimal policy for any point in time or space,  $u^*(t; x) = \arg \minou \maxod \big \langle \nabla_{\sti} V, f(x,u,d,t) \big\rangle$. %Clearly, this value holds great power, \textit{if} it can be solved.

\subsection{Hamilton-Jacobi-Bellman Solution for $V$}\label{sec:HJR}

Applying Bellman's principle of optimality to the value in \eqref{GameValue}, \cite{Evans1984DifferentialGames} proved that $V$ is equivalently the viscosity solution to the HJ-PDE
% \begin{Theorem}
% [Evans 84] \cite{evans1984differential} \\
% Given the assumptions (2.1)-(2.5) in [Evan 84], the value function $V$ defined in (\ref{GameValue}) is the viscosity solution to the following Hamilton-Jacobi Partial Differential Equation,
\begin{eqnarray}
\begin{aligned}
\dot{V}  + \min\{0, \:\ham(x, \nabla_{\sti} V, \tau)\} & = 0 &\text{ on } \ssp& \times \tint, \\
V(\sti, \tf) & = \tgfn(\sti) &\text{ on } \ssp &, 
\end{aligned}
\label{HJPDE-V}
\end{eqnarray}
where the Hamiltonian
% $\ham:\ssp \times \ssp \times [t,\tf] \mapsto \mathbb{R}$
is defined by
% \begin{eqnarray}
$\ham(x, p, t) = \minou \maxod \big\langle p, \sys(x, u, d, t) \big\rangle
\label{Hamiltonian}$.
% \end{eqnarray}
% \end{Theorem}

Thus, to solve the game, we may solve this HJ-PDE. In the years following this result, dynamic programming proved to yield high-fidelity solutions for a broad class of nonlinear systems and nonconvex games, making this a popular framework \textit{when the system dimension is low}. Due to the need for a discrete grid, dynamic programming methods suffer from the \textit{curse of dimensionality} and thus are infeasible for systems of dimension $n \geq 7$.

% In either game, to solve this PDE, therefore, yields the value function, the BRT and a robust optimal controller. The main challenge of HJR lies in solving this PDE in (\ref{HJPDE-V}); DP methods propagate $V(x, t)$ by finite-differences over a grid of points that grows exponentially with respect to $\statedim$ \cite{bansal2017hamilton}. In practice, this is computationally intractable for systems of $\statedim \geq 6$ and constrained to offline analysis.

% Additionally, the value functions can be used to derive the optimal control strategy for any point in space and time, e.g., for the \textit{Reach} game:
% \begin{eqnarray}
% \begin{aligned}
% u^*(t) = \arg \minou \nabla_{\sti} V(x,t) \cdot h_1(t)u.
% \end{aligned}
% \label{HJoc}
% \end{eqnarray}
% The main challenge of HJR lies in solving this PDE in (\ref{HJPDE-V}); DP methods propagate $V(x, t)$ by finite-differences over a grid of points that grows exponentially with respect to $\statedim$ \cite{bansal2017hamilton}. In practice, this is computationally intractable for systems of $\statedim \geq 6$ and constrained to offline analysis.

\subsection{Deep Learning of the HJ-PDE Solution}\label{subsec:hopf}

To avoid computation over a discrete grid, one may also use deep learning to approximate a PDE solution (\cite{raissi2017physics}). Namely, let the neural approximation $V_\theta$ be defined by
\begin{eqnarray}
\begin{aligned}
    V_\theta(\sti, \ti) &\triangleq \tgfn(\sti) + (\ti-\tf) \left(W_{Y} \big[\phi_{Y-1} \circ \phi_{Y-2} \circ ... \circ \phi_0 \big]\left([x,t]^\top \right) + b_{Y} \right), \\
    \phi_i(v) &\triangleq \sin (W_i v + b_i), \:\: i \in [1, Y]
\end{aligned}
\label{nnarc}
\end{eqnarray}
where $\phi_i:\mathbb{R}^{M_i} \to \mathbb{R}^{N_i}$ is the $i$-th layer of the neural network with weight matrix $W_i \in \mathbb{R}^{N_i \times M_i}$ and bias $b_i \in \mathbb{R}^{N_i}$ (\cite{sitzmann2020implicit}). To train this approximation, the PDE loss
\begin{eqnarray}
    \losspde(\theta; k) \triangleq \mathbb{E}_{\ssp, \tint_k} [\losstermpde(\theta)], \quad \losstermpde(\theta) \triangleq \Vert \dot{V}_\theta + H(x, \nabla_{\sti} V_\theta, t) \Vert,
    \label{def:PDEloss}
\end{eqnarray}
has proved successful, where $\tint_k \triangleq [\ti, \ti_k]$ is an increasing time range with $\ti_K\triangleq\tf$. This approach has proved to offer high-fidelity approximations of the solution $V$ for many systems beyond the dynamic programming limit (\cite{bansal2020deepreach}). However, in the limit of dimension, this method (like other learned PDE solutions) suffers from the complexity of optimizing the PDE loss, and the approximations can deteriorate.

\subsection{The Hopf Solution to HJ-PDE's with Linear Dynamics}\label{subsec:hopf}

Recently, \cite{darbon2016algorithms} recognized that when the dynamics in \eqref{prob} are defined by a linear function $\lsys(\mathsf{x},u,d,t) \triangleq A(\tau)\mathsf{x} + B_1 (\tau) u + B_2 (\tau)d$, an alternative to grid-methods is given by the Hopf formula, a solution to certain HJ-PDEs. Namely, if the game is given by
\begin{eqnarray}
\begin{aligned}
V_\lsys(\sti, \ti) &\triangleq \minot \supod \infou  \tgfn (\lintj(\tau)), \\
\end{aligned}
\label{LinGameValue}
\end{eqnarray}
where $\lintj$ is the trajectory of $\dot{\tj}_\lsys \triangleq \lsys(\lintj, \csig, \strd[u], \intt)$, $\tgfn$ is convex, and the linear Hamiltonian $H_\lsys(p,t) \triangleq \min_u \max_d \langle p, \Phi(\tf - t) (\nabla_u \lsys (\tf - t)u + \nabla_d \lsys(\tf - t) d) \rangle$ is convex (where $\Phi$ is the fundamental matrix of $\lsys$), then the value $V_\ell$ is given by 
\begin{eqnarray}
V_\ell (\sti,\ti) = -\min_{\tau \in \tint} \min_{p \in \mathbb{R}^{n}} \bigg\{ \tgfn^\star(p) - \langle \Phi(\tf - \tau)x, p \rangle + \int_0^{\tf-\tau} \ham_\ell (p, s) \:ds \bigg\} \triangleq \hopf[\tgfn, H_{\lsys}](x,t),
\label{HopfFormula}
\end{eqnarray}
where $J^\star(p):\ssp \rightarrow \mathbb{R} \cup \{+\infty\}$ is the convex-conjugate of $J$. This Hopf formula $\hopf$ may be \textit{independently} solved over space with non-smooth, convex optimization algorithms (proximal methods). This approach offers orders of magnitude of acceleration and memory efficiency over dynamic programming \cite{chow2017algorithm} by skirting the need for a grid, but is limited to the minimum over time problem (an upper-bound of $V$), and linear systems.
% Note, the minimum over time gives the union of Backwards Reachable Sets (BRS) \cite{bansal2017hamilton}, an under-approximation of the BRT.

\section{Approach}\label{sec:SfEnv}

We propose two methods for marrying the ability to solve the value $V_\lsys$ for simplified linear dynamics with the learned approach to approximating $V$ for the true nonlinear dynamics, and discuss the theoretical validity of these approaches.

\subsection{On the Relationship Between $V_\lsys$ and $V$}\label{subsec:inspo}

We begin by considering a few simple results which contextualize the linear form of the game. As is true in most dynamical systems theory, the game value for a nonlinear system is intimately related to the game value with a corresponding linearized version of the system. One perspective on this is given by the following results.

\begin{Theorem} \label{thm:conlinbd}
Let $\fset_c(t)$ be the $c$-level set of $V$ at $t$ and $\bar \fset(\tau)$ be a set containing any $\tj(s)$ s.t. $J(\tj(s')) \le c$ for $s,s' \in [\tau, \tf]$. Let
\begin{eqnarray}
%\begin{aligned}
    H^{\pm}_{\lsyse}(x,p,\tau) \triangleq H_\lsys(x,p,\tau) \pm \max_{\err \in \eset(\tau)} \pm \langle p, \err \rangle, %\\
    % % \quad \text{ and } \quad 
    % H^-_{\lsyse}(x,p,t) \triangleq H_\lsys(x,p,t) + \min_{\err \in \eset(t)} \langle p, \err \rangle,
    \label{errhams}
%\end{aligned}
\end{eqnarray}
where $\eset(\tau) \triangleq \{ \Vert \err \Vert \le \maxerr(\tau)\}$ and $\maxerr(\tau) = \max_\Sigma \Vert f - \ell \Vert$ is defined on $\Sigma \triangleq \bar \fset(\tau) \times \ctrlset \times \distset \times [\tau, \tf]$. 
If $V(x,\tau) = \min_{\tau'} \sup_{\mathfrak{d}} \inf_{\mathsf{u}} \tgfn (\tj(\tau'))$ for $x \in \fset_c(t)$ and $\tau \in \tint$, then
\begin{eqnarray}
% \begin{aligned}
    \big\vert V - V_{\lsys} \big\vert  \le \hopf[\tgfn, H_{\lsys+\err}^+] - \hopf[\tgfn, H_{\lsys+\err}^-] \triangleq \epsilon^*. \label{eq:conlinbd_eq}
    % H^+_{\lsyse} = H_\lsys + \max_{\err} \langle p, \err \rangle \:\: \& \:\: H^+_{\lsyse} = H_\lsys + \max_{\err} \langle p, \err \rangle.
% \end{aligned}
\end{eqnarray}
\end{Theorem}

%Although it is possible to rapidly solve these bounding linear game values and the bound itself, these values are irrelevant to this work. Rather, we consider this fact to 
See the Appendix in Sec.~\ref{apx} for the proof. Thm.~\ref{thm:conlinbd} offers a finite bound between the value of the linear and nonlinear games for any fixed time in a local region. 
% While this result applies to the minimum over time value, even when $V$ diverges from it, the bound in \eqref{eq:conlinbd_eq} is finite for any fixed time due to the boundedness of the values (\cite{mitchell2005time}). 
This knowledge justifies the use of the linear value in various scenarios, including its application in the current context for learning the nonlinear value. The following Corollary highlights an intuitive and desirable property of this bound: it diminishes in the limit of the operating point. 

\begin{Corollary} \label{cor:linbdtay}
With the linear game in (\ref{LinGameValue}) and assm. in Thm.~\ref{thm:conlinbd}, if $\lsys$ is defined by
\begin{eqnarray}
    \lsys(x,u,d,t) \triangleq  \nabla_xf|_{m_0}(t)x+ \nabla_uf|_{m_0}(t)u + \nabla_df|_{m_0}(t)d +\nabla_tf|_{m_0}t,
    \label{taylorlin}
\end{eqnarray}
with operating point $m_0 \triangleq (\tj_0, \csig_0, \dsig_0, \tau_0)$ and finite error $\delta^*(\tau) \sim O(\Vert m - m_0 \Vert^2)$, then
\begin{eqnarray}
    V(\sti,\ti) = V_\lsys(\sti,\ti) + \epsilon, \:\: \epsilon \in [0, \epsilon^*] 
    \quad \text{ s.t. } \quad 
    \lim_{m \to m_0} \epsilon^* = 0.
\end{eqnarray}
\end{Corollary}

Thm.~\ref{thm:conlinbd} and Cor.~\ref{cor:linbdtay} imply that computing $V_\lsys$ provides a portion of the solution that is dominant near the operating point. 
% Despite significant differences with $V$, $V_\lsys$ nonetheless offers accurate reinforcement for the early time window, a period known to be violated by local optima in the learning program.
For these reasons, we propose using the linear solution in the nonlinear program, specifically with a least-squares supervision loss.

% \subsection{Supervision with $V_l$}

% To leverage the linear solution, we define the following linear supervision loss.
\begin{Definition} [Linear Supervision Loss] Let the linear supervision loss be defined by
    \begin{eqnarray}
    \losstermls(\theta) \triangleq \rho \Vert V_\theta - V_\lsys \Vert + \rho_g\Vert \nabla_x V_\theta - \nabla_x V_\lsys \Vert
    \end{eqnarray}
    where $\rho, \rho_g \in \mathbb{R}$ are hyperparameters.
    \label{def:LSloss}
\end{Definition}

The least-squares loss satisfies the local Polyak-{\L}ojasiewicz (PL) condition for wide models with smooth activations (\cite{liu2022loss}).
This property guarantees exponential convergence to a global optimum with (stochastic) gradient descent, an increasingly beneficial fact as the dimension rises.

\subsection{Decayed Linear Semi-Supervision for Learning $V$}\label{subsec:LSSd}

In light of the relation between $V$ \& $V_\lsys$ in Cor.~\ref{cor:linbdtay}, and the simplicity of learning $V_\lsys$ with supervision, we propose the following simple semi-supervision loss, combining linear supervision with the PDE loss term to learn the nonlinear value. 

\begin{Definition} [Linear Semi-Supervision Loss, Decayed] Let a loss be defined by
    \begin{eqnarray}
    \loss_{LSS-D}(\theta; k) \triangleq (1 - \augp_k) \mathbb{E}_{\ssp, \tint}[\losstermls(\theta)] + \augp_k \mathbb{E}_{\ssp, \tint} [\losstermpde(\theta)], \quad \text{ where } \quad\augp_0\triangleq0.
    \end{eqnarray}
    Here $k$ is the iteration and $\augp_k$ is a monotonically increasing value starting at $\augp_0\triangleq0$.
    \label{def:LSSDloss}
\end{Definition}

The key idea is to use $\lambda_k$ to gradually transition the neural value from approximating the linear solution $V_\ell$, a global minimizer of $\losstermls$, to approximating the nonlinear solution $V$, a global minimizer of $\losstermpde$.
We note there is no guarantee of a globally minimizing path between $\losstermls$ and $\losstermpde$ for $\augp_0 \to \augp_K$, but empirically this works well nonetheless, as shown in Sec.~\ref{sec:demos}. Rather than learn the solution from scratch, one might note that the program only needs to correct a partially true solution and, thus, has a simpler task at hand. This also allows one to avoid the gradual time-curriculum approach that previous authors have found necessary for learned PDE solutions (\cite{bansal2020deepreach}), yielding a highly accelerated and performant program.

% , and yet, in experimentation, this program boasts unrivaled speed and a significant improvement in accuracy over the baseline approach in high-dimensions, perhaps because of the ``niceness'' of $\losstermls$. 
% In fact, while $V$ is a global minimizer of $\augp_K=1$, we found that $\augp_K<1$ offers much better accuracy in higher-dimensions, suggesting the importance of the supervision term to structure the learning problem beyond solely preconditioning.

\subsection{Linear Semi-Supervision in an Augmented Game for Learning $V$}\label{subsec:LSSaug}

\begin{figure}[t]
    \centering
    \includegraphics[width=\linewidth]{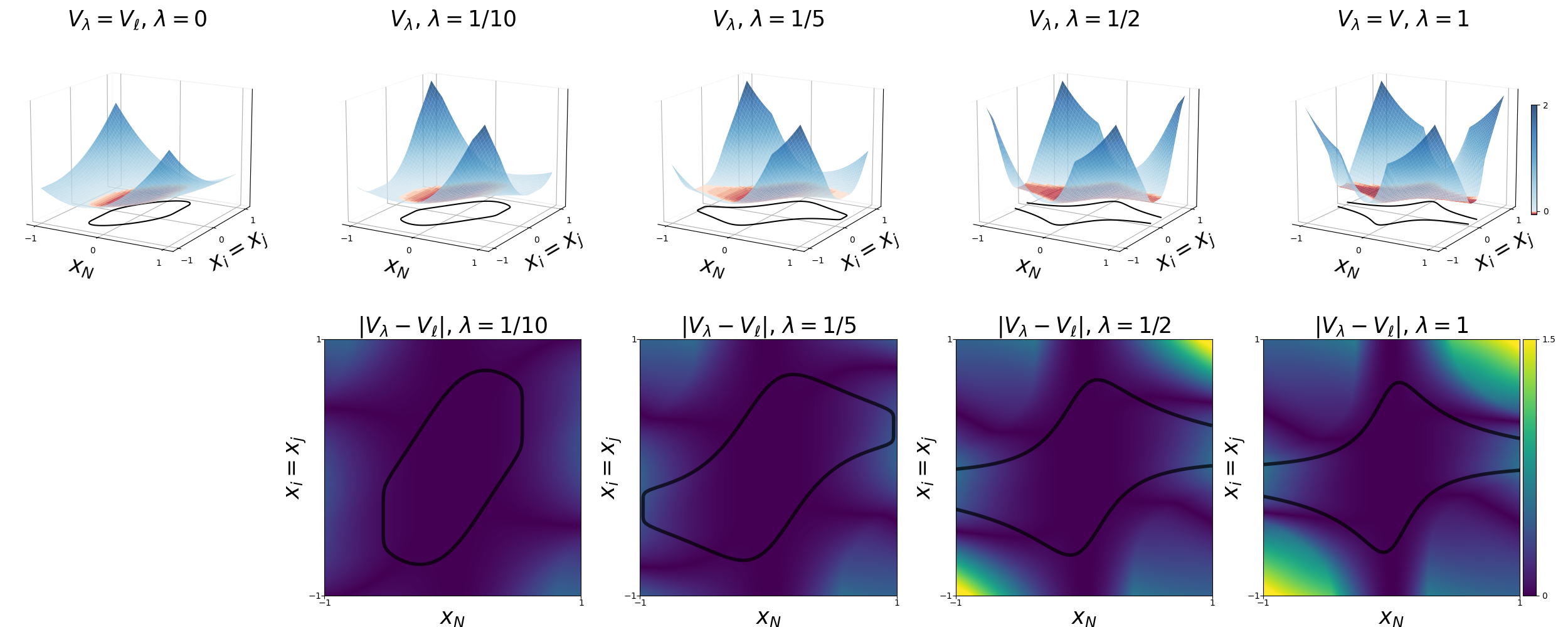}
    % {figs/Lambda_Variation_Demo_Plot_v1_1_viridis.png}
    \caption{\textbf{Demonstration of Thm.~\ref{thm:Vlam} [$V_\augp$] and Cor.~\ref{cor:linbdtay}} On top, the true value of $V_\augp$ at $t=1$ for the problem posed in \eqref{def:auggame} with $N=3$ along the range of $\augp$ is given. Note the smooth change from $\augp=0$, where $V_\augp = V_\lsys$, to $\augp=1$, where $V_\augp = V$. In the bottom row, the error between $V_\augp$ and $V_\lsys$ is plotted as $\augp$ increases. Note the gradual increase in error and the large regions of $V$ with low error.}
    \label{fig:theorydemo}\vspace{-1em}
\end{figure}

In addition to the decayed LSS program, we propose a more sophisticated method that learns the solution of an augmentation of the game. By defining a continuous spectrum of systems between the nonlinear and linear systems, we show that the augmented game offers a setting where $\losstermls$ and $\losstermpde$ are not in conflict. 
% In the benchmark experiments, this program succeeds on the timescale of the baseline program, and hence is slower than the decay program from Sec.~\ref{subsec:LSSd} (likely due to the increased complexity), but it \textit{yields the most accurate solutions}. 
Consider the following augmented dynamics.
\begin{Definition} {\bf{Nonlinear Spectrum Augmentation}} \\
    Let an augmented system be defined by state $\augtj = [\tj, \lambda]^\top$ in $\sasp \triangleq \ssp \times \mathbb{R}$ with dynamics
    \begin{eqnarray}
        \dot{\augtj} \triangleq \begin{bmatrix}
            (1 - \augp) \lsys (\tj, u, d, t) + \augp\sys (\tj, u, d, t) \\ 0
        \end{bmatrix} \triangleq \augsys( \augtj, u,d,t).
        \label{def:augsys}
    \end{eqnarray}
    Let the augmented game be defined with the same inputs and strategies as $V$, s.t.
    \begin{eqnarray}
        V_\augp(x, \lambda, t) \triangleq \supod \infou \minot \tgfn (P\augtj(\tau)).
        \label{def:auggame}
    \end{eqnarray}
    where $P:\sasp \to \ssp$ represents projection from the augmented space given by $P([x, \lambda]^\top)= x$.
\end{Definition}
Intuitively, this system represents a continuous spectrum of nonlinear systems between any given linear and nonlinear dynamics. Notably, this game has a few simple properties that will make it valuable for learning.
 \begin{Theorem} \label{thm:Vlam}
$V_\lambda$ is the viscosity solution of
\begin{eqnarray}
\begin{aligned}
    \dot{V}_\lambda  + \min\{0, \ham_\lambda(\tilde{\sti}, \nabla_{\sti} V_\lambda, t)\} & = 0 &\text{ on } \sasp& \times \tint, \\
    V_\lambda(\tilde{x}, \tf) & = \tgfn(\sti) &\text{ on } \sasp,&
\end{aligned}
\label{HJPDE-Vlam}
\end{eqnarray}
where $H_\augp(x, \augp, p, t) \triangleq \minou \maxod \big\langle p, (1 - \augp)\ell(x,u,d,t) + \augp f(x,u,d,t) \big\rangle$.
% \begin{eqnarray}
% \begin{aligned}
%     H_\augp(x, \augp, p, t) = \minou \maxod \big\langle p, (1 - \augp)\ell(x,u,d,t) + \augp f(x,u,d,t) \big\rangle.
% \end{aligned}
% \label{Haug}
% \end{eqnarray}
Moreover, if $V_\ell = \min_{\tau} \sup_{\mathfrak{d}} \sup_{\mathsf{u}} \tgfn (\lintj(\tau))$,
\begin{eqnarray}
\begin{aligned}
    V_\augp(x, \lambda, t) = \begin{cases} V_\ell(x,t)  & \text { if } \: \augp = 0, \\
    V(x,t)  & \text { if } \:\augp = 1, \\
    \end{cases}
\end{aligned}
\label{auggamebc}
\end{eqnarray}
and $V_\augp$ Lipschitz continuous w.r.t. $\augp$.
\end{Theorem}
We provide the formal proof in the Appendix (Sec.~\ref{apx}) and a demonstration in Fig.~\ref{fig:theorydemo}. Due to the invariance of the $\augp$ trajectory with respect to the inputs, the game is independent for different values of $\augp$. By this, we mean that the players cannot affect $\augp$, and hence the optimal strategies are preserved on $\augp$ slices. Yet there is a continuity in the game value across $\augp$ that will be valuable for simplifying the learning problem. 

The purpose of the augmented game is to host a program in which the linear supervision and PDE losses have an intersecting set of global minimizers. The work of \cite{evans1984differential} certifies that $V_\augp$ is the solution of the HJ-PDE in (\ref{HJPDE-Vlam}), and hence may be learned with $\losstermpde$. With (\ref{auggamebc}), we propose limiting the supervision loss to the $\augp=0$ subspace where $V_\augp= V_\ell$, while applying the PDE loss to the entire augmented state space to solve a continuous solution across $\augp$. This is given explicitly by the following loss.

\begin{Definition} [Linear Semi-Supervision Loss, Nonlinear Spectrum] Let a loss be
    \begin{eqnarray}
    \loss_{LSS-NS}(\theta; k) \triangleq \mathbb{E}_{\tilde{\ssp}_{\augp=0}, \tint_k}[\losstermls(\theta)] + \mathbb{E}_{\tilde{\ssp}, \tint_k} [\losstermpde(\theta)]
    \end{eqnarray}
    where $\tint_k \triangleq [\ti, \ti_k]$ and $\ti_k$ is a monotonically increasing time range with $\ti_0\triangleq\ti$ and $\ti_K\triangleq\tf$.
    \label{def:LSSNCloss}
\end{Definition}

Note, in this case, $V_\theta$ has an input of one increased dimension to accommodate $\augp$. In some sense, this program adds structure to the problem by incorporating another boundary condition (corresponding to $V_\lsys$ at $\augp = 0$) that allows us to ultimately better approximate $V$, which lives on the $\augp=1$ subspace of $V_\augp$.% Accordingly, it might be possible to learn solely the difference between $V$ and $V_\ell$.

\subsection{Generating the Linear Supervisor $V_\lsys$}\label{subsec:genVL}

%When a system is linear, the strength of the Hopf formula is that it offers a rapid method to solve $V_\lsys$ and the corresponding optimal control \textit{at an isolated point}, making it far more efficient than dynamic programming. To learn a function with input in $\mathbb{R}^{n+1}$, however, the cost of computation is on the order of milliseconds, a little slow for sampling large batches of the space continuously. Additionally, the Hopf formula is itself solved by proximal algorithms that may yield noisy data with insufficient parameters. 
%We found that a fast and accurate method is to solve a finite bank of Hopf data, train a neural net on this linear data via supervision, and then introduce the PDE loss after some period to smooth the supervised approximation. While slightly slower, we also found that training a neural net with the PDE loss alone offers a fairly accurate linear approximation, which is useful if the Hopf data generation is infeasible. In either case, the learned model for $V_\lsys$ is fixed afterward and then queried rapidly across the state space in the following programs to learn $V$.  \wsnote{PARE}

Here we introduce two options for constructing the linear supervisor $V_\lsys$. First, the Hopf formula \eqref{LinGameValue} can generate samples of the linear value, which are then used to train a neural network via supervision, employing the PDE loss $\losstermpde$  \eqref{def:PDEloss} to encourage smooth solutions. Alternatively, one can train the neural network directly with $\losstermpde$ without the Hopf formula. This approach tends to be slower, but similar in accuracy.  The learned model $V_\lsys$ is then fixed and queried across the state space in the following programs to learn $V$. 

\section{Demonstration}\label{sec:demos}

\subsection{$N$-Dimensional, Differential Game Benchmark}\label{sec:benchproblem}

% \begin{figure}[t]
%     \centering
%     \includegraphics[width=\linewidth]{fig_sketches/BM_demo_figsk.jpg}    
%     \caption{\textbf{$N$-d Benchmark Demo} }
%     \label{fig:TVonly}\vspace{-1em}
% \end{figure}

% \begin{figure}[t]
%     \centering
%     \includegraphics[width=\linewidth]{figs/Lambda_Variation_Demo_Plot_final_nozls_lincolor.png}
%     % {figs/Lambda_Variation_Demo_Plot_v1_1_viridis.png}
%     \caption{\textbf{Demo of Error and $V_\augp$ Propositions in the $N$-d Benchmark} Here $N=3$, $t=1$ and at $\augp=1$ $(\alpha, \beta, \gamma) = (-20,-20,1)$. The top row shows how the value transforms from the linear (left) to the nonlinear solution (right) in any case. Below is the difference between the linear and nonlinear values ($\epsilon$) for each value of $\augp$.}
%     \label{fig:TVonly}\vspace{-1em}
% \end{figure}

% \begin{RunningExample} ($N$-dimensional, Nonlinear Publisher-Subscriber System) \\

For demonstration, we first introduce a ``publisher-subscriber'' game, chosen such that the ground truth solution can be directly computed (via dynamic programming) for comparison. Consider a system with a ``publisher'' state $x_0$ which unidirectionally influences ``subscriber'' states $x_i$, such that
\begin{eqnarray}
    \begin{bmatrix} \dot \tj_0 \\ \dot \tj_i \end{bmatrix} \triangleq \begin{bmatrix} a & 0 \\ -1 & a \end{bmatrix} \begin{bmatrix}  \tj_0 \\  \tj_i \end{bmatrix} + \begin{bmatrix} 0 \\ b \end{bmatrix}u_i + \begin{bmatrix} 0 \\ c \end{bmatrix}d_i + \begin{bmatrix} \alpha \sin(\tj_0)\tj_0^2 \\ -\beta  \tj_0 \tj_i^2 \end{bmatrix},
\label{pubsub2d}
\end{eqnarray}
with $u_i \in \{ \vert u_i \vert \le 1 \}$, $d_i \in \{\vert d_i \vert \le 1 \}$ and $a, b, c, \alpha, \beta \in \mathbb{R}$. 
In this game, let the agent's goal be to attenuate the subscribers while the opponent will seek to amplify them. 
Hence, let the target set be a ball of radius $r$ with $J_{\mathcal{T},i}(x_0, x_i) \triangleq \frac{1}{2}(x_0^2 + x_i^2 - r^2)$ and let the goals of the agent and the opponent be to minimize and maximize $J_{\mathcal{T},i}$ respectively.
\begin{figure}[t]
    \centering
    \includegraphics[width=\linewidth]{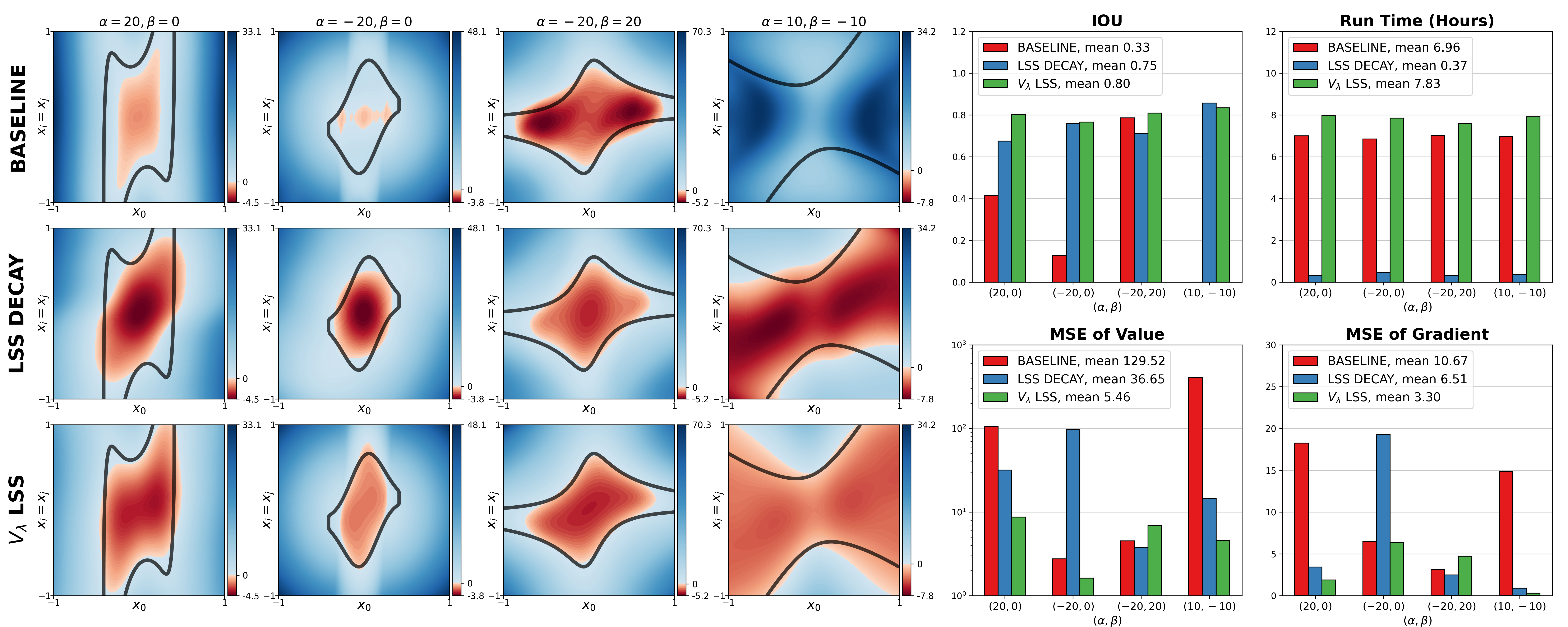}    
    \caption{\textbf{50-D Benchmark Result Comparison} On the left, a slice of the learned solution for four variations (columns) of \eqref{pubsubNd} where $(\alpha, \beta) \in \{(20, 0), (-20, 0), (-20, 20), (10, -10)\}$ is shown for each proposed method (rows), and the ground truth zero-level set is overlaid in black. On the right, the IOU, MSEs, and run time are given for each of the variations and methods.
    }
    \label{fig:benchresults}\vspace{-1em}
\end{figure}

With $N-1$ subscribers, the composed system takes the form 
\begin{eqnarray}
    \dot{\tj} = \big(e_1 e_1^\top - \mathbf{1}_{N} e_1^\top + a I_N \big) \tj + \begin{bmatrix} \mathbf{0}_{N-1}^\top  \\ bI_{N-1} \end{bmatrix}u + \begin{bmatrix} \mathbf{0}_{N-1}^\top \\ cI_{N-1} \end{bmatrix}d + \begin{bmatrix} \alpha \sin(\tj_0) \\ -\beta \tj_0 \mathbf{1}_{N-1} \end{bmatrix} \circ (\tj \circ \tj),
\label{pubsubNd}
\end{eqnarray}
where $\tj = [ \tj_0, ..., \tj_{N-1}]^\top \in \mathbb{R}^N$, $u \in \{ u \in \mathbb{R}^{N-1} \mid \Vert u \Vert_\infty \le 1 \}$, $d \in \{d \in \mathbb{R}^{N-1} \mid \Vert d \Vert_\infty \le 1 \}$, and here $\circ$ represents element-wise multiplication. 
Moreover, the composed target may be written as $\tgfn(x) = \frac{1}{2} \left((N-1)x_0^2 + \sum x_i^2 - (N-1)r^2 \right) = \sum J_{\mathcal{T},i} (P_i x)$ where $P_i:\mathbb{R}^N\to\mathbb{R}^2$ is projection from the composed space to each $2$-D publisher-subscriber space.
\begin{Remark} \label{rem:benchprop}
    The value of this game has the special properties,
    \begin{eqnarray}
    \begin{aligned}
    V(x, t) = \sum^{N-1} V_i(P_i x, t), \quad \& \quad  \hjrset(t) \cap \tilde {\mathcal{X}} = \hjrset_i(t),
    \end{aligned}
    \label{pubsubval}
    \end{eqnarray}
    where $\tilde {\mathcal{X}} \triangleq \{\tilde x \in \mathbb{R}^N \mid \forall i,j > 0, \tilde{x}_i = \tilde{x}_j \}$ is a diagonal. For proof, see the Appendix (Sec.~\ref{apx}).
\end{Remark}
Hence, we may solve the value of the $N$-D game value by summing the values of the $N-1$ decomposed 2-D games, which are solved with dynamic-programming. This is beneficial for a high-dimensional assay as we may naively solve the game in the full system (\ref{pubsubNd}) with a given method and score it on the high-fidelity dynamic programming solution.

To interrogate the proposed methods, let $N=50$ and consider four nonlinear parameter variations of \eqref{pubsubNd}. For each, we train a 3-layer $V_\theta$ with the baseline program ($\losspde$, \eqref{def:PDEloss}), the decayed linear semi-supervision program ($\losslssd$, Def.~\eqref{def:LSSDloss}), and the augmented game $V_\augp$ program ($\losslssnl$, Def.~\eqref{def:LSSNCloss}). See the Appendix in Sec.~\ref{apx} for training details. We score them against the dynamic programming solution with the Intersection-Over-Union (IOU) ($|\hjrset_1 \cap \hjrset_2|/|\hjrset_1 \cup \hjrset_2|$) and the Mean-Squared Error (MSE) of both the value and the value gradient. The results are displayed in Fig.~\ref{fig:benchresults}.

% Furthermore, both the nonlinearity and control authority of this model can be tuned by choosing various values of parameters, e.g., with $\alpha, \gamma=0$ the problem is linear and with $c=0$ the differential game reduces to an optimal control problem. \wsnote{PARE}

% Lastly, we note that in this setting the ground truth also offers a method of performing a ``capacity'' test of the neural net, in which we may solve the maximum ability of a given architecture to approximate the true solution by direct supervision,
% \begin{eqnarray}
% \loss_{cpy}(\theta) = \mathbb{E}_{\ssp,\tint}[\Vert V_\theta - V_{gt} \Vert]
% \end{eqnarray}
% where $V_{gt}$ is computed from (\ref{pubsubval}) by summation of decomposed dynamic-programming solutions. Given that this capacity purely uses the supervision loss, we expect near-optimal approximation of the converged program.
The results demonstrate that both proposed approaches offer significant improvement over the baseline. The augmented program with $V_\augp$ achieves the highest mean IOU (2.4x baseline) and lowest mean MSEs (23.7x baseline) with a slight cost in run time over the baseline (1.1x baseline). Alternatively, the decayed program yields fairly accurate IOU (2.2x baseline) and MSEs (3.5x baseline) with a nearly 20-fold acceleration over the baseline, since no time curriculum is required during training. Ultimately, both programs underscore the benefit of the supervision loss in learning the HJ-PDE solution in high dimensions. 

\subsection{10D Quadrotor Optimal Control for Collision Avoidance}\label{sec:benchproblem}

\begin{table}[b!]
\caption{10D Quadrotor Results}
\centering
\begin{tabular}{|p{6cm}|p{3cm}|p{1.5cm}|p{1.5cm}|p{1cm}|}
\hline
Method & Recovered Volume & FP\% & FN\% & Time \\
\hline
Baseline (\cite{bansal2020deepreach}) & 81.09\% & 1.864\% & 0.005\% & 2.5h \\
LSS Decay (Def.~\ref{def:LSSDloss}) & 94.89\% & 0.725\% & 0.049\% & 0.5h \\
LSS Decay Adp. (App. \ref{appendix:adaptive_weight}) & 94.96\% & 0.234\% & 0.326\% & 0.8h \\
$V_\lambda$ LSS (Def.~\ref{def:LSSNCloss}) & 88.92\% & 0.908\% & 0.041\% & 4.5h \\
\hline
\end{tabular}
\label{table:quadrotor_results}
\end{table}

Let a drone be flying with high-velocity toward an obstacle that it would like to avoid. Consider a 10-D quadrotor dynamic model (\cite{gong2024robustcontrollyapunovvaluefunctions}) with state defined as $x = [p_x, v_x, \theta, \omega_y, p_y, v_y, \phi, \omega_x, p_z, v_z]^T$. Here, $(p_x, p_y, p_z)$ represent the position, $(v_x, v_y, v_z)$ are the velocities in the world frame, ($\phi$, $\theta$) denote the roll and pitch angles, and $(\omega_x, \omega_y)$ are the corresponding angular rates.
The full dynamics can be found in the Appendix (Sec.~ \ref{appendix:quad_dyn}).
The target set is given by $\mathcal{T}=\{ x \mid p_x^2+p_y^2 \leq 0.5^2 \}$, representing a cylindrical obstacle.

The drone dynamics are linearized around the trivial operating point and then used to learn the linear value function $V_\ell$ with only the PDE loss $\losspde$. 
With $V_\ell$, we learn the nonlinear solutions through both proposed semi-supervised approaches. We compare these to the baseline method proposed in \cite{bansal2020deepreach}, which solely uses $\losspde$. 
As tuning \eqref{def:LSloss} can be time-consuming, we also introduce an adaptive-weighting scheme (see Appendix, Sec.~\ref{appendix:adaptive_weight}) that may work for either approach but is applied here to the decay method.
The performance gain is quantified with three metrics: the volume of the corrected safe set through probabilistic conformal expansion (referred to as recovered volume), roll-out false-negatives (FN\%) of the solution (without conformal expansion), and roll-out false-positives (FP\%). Intuitively, $FP\%$ measures the proportion of overly optimistic classifications while $FN\%$ measures the proportion of overly conservative classifications. Detailed descriptions of these metrics are provided in the Appendix (Sec.~\ref{appendix:metrics}).

\begin{figure}[t]
    \centering
    \includegraphics[width=0.85\linewidth]{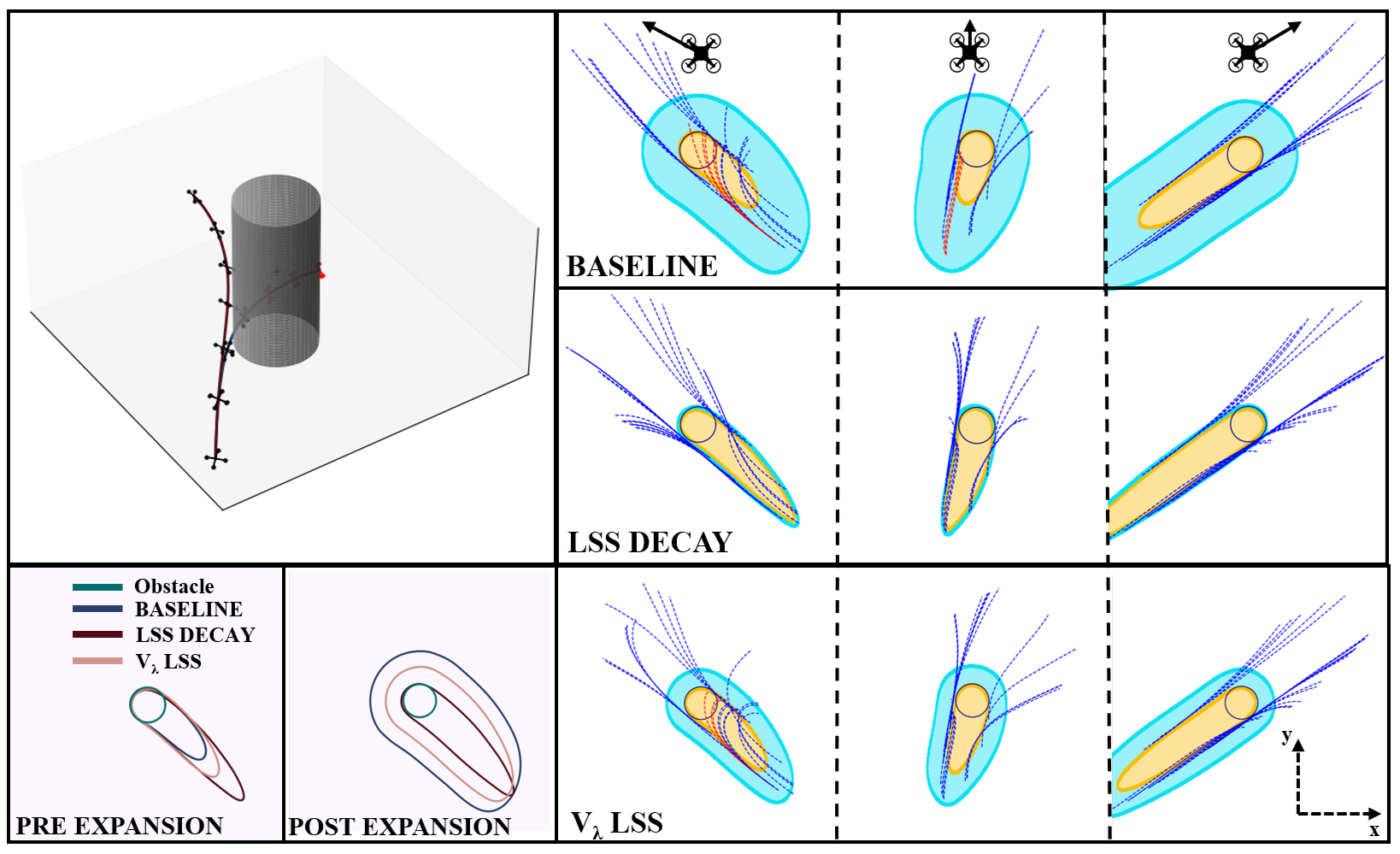}
    \caption{\textbf{10-D Quadrotor Result Comparison} In the upper-left, the problem in which the drone is flying toward an obstacle is depicted with two trajectories demonstrating success and failure. On the right, slices of the sub-zero level set of the learned value that approximate the unsafe set are shown (gold), along with the 99.9\%-confidence conformal expansion of the learned set (teal) and a sample of the roll-outs (blue if safe, else red). In the lower-left, a slice of the learned sets is shown before and after the conformal expansion.}
    \label{fig:quadrotor_results}
\end{figure}

The results are summarized in Fig.~\ref{fig:quadrotor_results} and Table~\ref{table:quadrotor_results}. The proposed variants significantly outperform the baseline model in terms of recovered volume and false positives. All the methods maintain a reasonable false-negative rate, as FN\% is less critical from a safety perspective --- false negatives correspond to issuing false alarms, whereas false positives can lead to system failures. 

With all methods, FP\% is significantly higher than FN\%, corresponding to an under-approximation of the unsafe set, perhaps due to a bias from the high-value in high-dimensions. It appears linear semi-supervision protects against this bias, grounding the program in some sense. Interestingly, the faster LSS decay methods demonstrate superior performance compared to $V_\lambda$ LSS. We hypothesize that this is because the true solution is not particularly nonlinear. As a result, the decay schemes perform effectively by ``polishing'' the linear solution, while the augmented system introduces additional complexity. In scenarios where the nonlinear solution deviates significantly from $V_\lsys$, however, the augmented variant may be more robust. \vspace{-1em}
% In reality, as discussed in Sec.~\ref{subsec:genVL}, this linear solution is generated from an optimization process itself and thus can be noisy. 

% \begin{figure}[t]
%     \centering
%     \includegraphics[width=\linewidth]{fig_sketches/QR_RO_LS_figsk.jpg}
%     \caption{\textbf{Quadrotor Results} WE HAVE SOME GOOD EXPANSION PLOTS FOR THESE, ZEYUAN MAKING A 3D PLOT \textcolor{mygray}{\lipsum[4]}}
%     \label{fig:TVonly}\vspace{-1em}
% \end{figure}

% \begin{table}
% \caption{10D quadrotor results}
% \begin{tabular}{ |p{3.5cm}|p{3.1cm}|p{2.3cm}|p{2.35cm}|p{1cm}|  }
% \label{table:quadrotor_results}
% \hline
% % Method & recovered volume &FP\% & FN\% & time \\
% \multicolumn{5}{|c|}{Method & recovered volume & FP\% & FN\% & Time} \\
% \hline
% Baseline & 81.09\% &1.864\% & 0.005\% &2.5h\\
% LSS Decay & 94.89\% & 0.725\% & 0.049\%& 0.5h \\
% LSS Adaptive & 94.96\% & 0.234\% & 0.326\% & 0.8h\\
% $V_\lambda$ LSS & 88.92\% & 0.908\% & 0.041\% & 4.5h\\
% \hline
% \end{tabular}
% \end{table}

\section{Conclusion}\label{sec:conclusion}

In this work, we propose and demonstrate the introduction of the HJ-PDE \textit{linear} solution to the problem of learning high-dimensional, \textit{nonlinear} HJ-PDE solutions for control and differential games. This proves to be a significant insight, offering improvements in speed in a decayed method and accuracy in an augmented method. For a practitioner, the experimental results suggest that if the problem is very nonlinear, the augmented method performs best, while the decay suffices otherwise (and is greatly accelerated). We note the framework supports any linearization and we offer theory about the Taylor linearization specifically. Moreover, it seems one might be able to improve our results with an ensemble of linear models but we leave this to future work. We believe these results are a valuable step for high-dimensional approaches to autonomy, and plan to extend this to real world applications.

\section*{Aknowledgements}\label{sec:aknowlegements}

This work was supported by ONR grant N00014-24-1-2661 and NSF CAREER (Award \# 2240163). We would like to thank Albert Lin, Dylan Hirsch, Feng Yi, Nikhil Shinde, Sander Tonkens, and Yat Tin Chow for helpful discussions related to the paper.

\bibliography{l4dc_2024/references, l4dc_2024/refs_herbert}
% \bibliography{main}

\section{Appendix} \label{apx}

\subsection{Proof of Theorem~\ref{thm:conlinbd}}

\begin{proof}
    First note, by the triangle inequality,
    $$\big\vert V(x,t) - V_{\lsys}(x,t) \big\vert \le \big\vert \bar V(x,t) - V_{\lsys}(x,t) \big\vert + \big\vert \bar V(x,t) - V(x,t) \big\vert,$$
    where $\bar V(x,t) \triangleq \minot \supod \infou J(\tj(\tau))$ is the minimum over time value, an upperbound of $V$. By assumption, we consider $(x,t)$ s.t. $|V(x,t) - \bar V(x,t)| = 0$, which corresponds to states where the optimal disturbance strategy doesn't vary with the time of arrival. Beyond this assumption, the right-hand bound in finite for any $\ti>-\infty$ and approaches $0$ as $\ti \to \tf$, since both values are bounded (\cite{mitchell2005time, evans1984differential}). It remains to prove the boundedness of $\big\vert \bar V(x,t) - V_{\lsys}(x,t) \big\vert = \big\vert V(x,t) - V_{\lsys}(x,t) \big\vert$.
    
    To prove this, we first generalize the Theorem 1 in \cite{sharpless2024conservative} to include a time-varying endpoint for the minimum over time value. We also show this holds in an error-assisting game, which gives the lower bound. We solely consider the reach game for brevity but as in \cite{sharpless2024conservative}, mirrored results hold in the avoid game.

    By assumption, $\bar \fset(\tau)$ is any set s.t. for any user-defined $c \in \mathbb{R}$,
    \begin{eqnarray}
    % \begin{aligned}
    \bar\fset(\tau) \supseteq \{ y \in \ssp  \mid  \:\: y = \tj(s; x, \csig, \dsig, t), \text{ s.t. } \tgfn (\tj(s')) \le c,\: s,s' \in [\tau, \tf]\},
    \label{errenvhypo}
    % \end{aligned}
    \end{eqnarray}
    meaning $\bar\fset(\tau)$ covers all trajectories which cross the $c$ level of the target up to backwards time $\tau$. The maximum error is then given by $\maxerr(\tau) \triangleq \max_{\bar \fset(\tau) \times \ctrlsetsig \times \distsetsig \times \tint} \Vert [f - \lsys] (x, u, d, \tau) \Vert$. 

    Consider the linear affine dynamical system, 
    \begin{eqnarray}
        \dot \tj_{\lsyse} (\tau) = \ell(\tj(\tau), \csig(\tau), \dsig(\tau), \tau) + \esig (\tau),
    \end{eqnarray}
    where $\esig : \tint \to \eset $ is the measurable signal of an error \textit{player} with actions in $\eset \triangleq \{ \err \in \ssp \mid \Vert \err \Vert \le \maxerr(\tau) \}$. Moreover, let this player have non-anticipative strategies $\stre: \ctrlsetsig \to \esigset$ that map from control signals to error signals.
    
    Let $V_{\lsyse}^+$ and $V_{\lsyse}^-$ be the values of the games over the linear dynamics with error where the error player antagonizes or assists the controller respectively such that
    \begin{eqnarray}
        \begin{aligned}
        V_{\lsyse}^+(x,t) &\triangleq \minot \supoe \supod \infou \tgfn ( \tj_{\lsyse}(\tau) ), \\
        V_{\lsyse}^-(x,t) &\triangleq \minot \infoe \supod \infou \tgfn ( \tj_{\lsyse}(\tau) ).
        \end{aligned}
    \end{eqnarray}
    In the zero-sum setting, the error player aligns wth the disturbance or the control respectively and may equivalently be considered as an augmentation of their potentials.
    
    Consider the following facts.
        
    \begin{Lemma} For any $x \in \fset_c (t)$,
    $V(x,t) \le V_{\lsyse}^+(x,t)$.
    \label{lem:conservupper}
    \end{Lemma}
    \begin{proof}
        The proof is given by the construction a specific error strategy that yields the nonlinear trajectory which the error player may induce to its benefit.

        We aim to show that for any $t$, $x \in \fset_c(t)$, and $c' \le c$,
        \begin{eqnarray}
            V(x, t) > c' \implies V_{\lsyse}^+(x, t) > c'.
        \end{eqnarray}
        By definition $x \in \fset_c(t) \implies V(x,t) \le c$, hence, it follows $x \in \fset_c(t) \implies V(x,t) \le V_{\lsyse}^+(x,t)$.

        Assume 
        \begin{eqnarray}
            c' < V(x, t) = \minot \supod \infou \tgfn (\tj( \tau)).
        \end{eqnarray}
        Since $V$ is continuous w.r.t. $\tau$ (\cite{Evans84}), the infimum is attained, and thus $\forall \tau \in \tint$,
        $$ c' < \supod \infou \tgfn (\tj (\tau; \csig, \strd[\csig], t)).$$
        $\implies  \exists \sigma > 0, \exists \tilde \strd \in \stratset$ yielding $\tilde \tj$
        \begin{equation}
            c' < c' + 2\sigma < \infou \tgfn (\tilde \tj (\tau; \csig, \tilde\strd[\csig], t)), 
        \end{equation}
        $\implies  \exists \sigma > 0, \exists \tilde \strd \in \stratset, \forall \csig \in \ctrlsetsig$
        \begin{equation}
            c' < c' + \sigma < \tgfn (\tilde \tj (\tau; \csig, \tilde\strd[\csig], t)) \le c,
            \label{inJc}
        \end{equation}
        where the upper bound is assumed for $\tilde \strd$ w.l.o.g. as $V(x,t) < c$. 
        
        Let the ``true'' error signal $\hat {\esig}$ for $\tj$, $\csig$ and $\dsig$ be defined as
        \begin{eqnarray}
            \hat {\esig} (\tau) \triangleq \sys(\tj(\tau), \csig(\tau), \dsig(\tau), \tau) - \lsys(\tj(\tau), \csig(\tau), \dsig(\tau), \tau).
        \end{eqnarray}
        Moreover, for any $\csig$, let the non-anticipative $\hat \stre$ be s.t. $\hat \stre [\csig] = \hat \esig$, inducing the ``true'' trajectory $ \tj$ corresponding to $\csig$. Note, for $\tilde \strd$, the true trajectory $\tilde \tj$ satisfies $$\tgfn(\tilde \tj(\tau)) \le c \implies \hat \esig(\tau) \in \bar \fset(\tau),$$ 
        therefore, it follows from Theorem 3 in \cite{sharpless2024conservative} that $\exists \tilde \tj_{\lsyse}$ s.t. $$\exists \tilde \tj_{\lsyse}(\tau; \csig, \tilde \strd[u], \hat \stre [u], t) = \tilde \tj (\tau; \csig, \tilde \strd[u], t), \:\forall \tau \in \tint.$$
        Then it must be that
        \begin{equation}
            c' < \infou \tgfn (\tilde \tj_{\lsyse}(\tau; \csig, \tilde \strd[u], \hat \stre [u], t)) < \supoe \supod \infou \tgfn (\tj_{\lsyse}(\tau; \csig, \strd[u], \stre [u], t)), 
        \end{equation}
        and since this holds for any time $\tau$,
        \begin{equation}
            c' < \minot \supoe \supod \infou \tgfn (\tj_{\lsyse}(\tau; \csig, \strd[u], \stre [u], t)) = V_{\lsyse}^+(x,t). 
        \end{equation}
    \end{proof}

    \begin{Lemma} For any $x \in \fset_c (t)$,
    $V(x,t) \ge V_{\lsyse}^-(x,t)$.
    \end{Lemma}
    \begin{proof}
        By Isaac's condition one observes
        $$V(x,t) = \minot \supou \infod \tgfn ( \tj(\tau) ).$$
        The proof is then analogous to that of Lemma~\ref{lem:conservupper} by constructing $\hat\esig$ and $\hat \stre$ for a $\tilde \stru$ and any $\dsig$.
    \end{proof}

    \begin{Lemma} For any $x \in \ssp $,
    % $V(x,t) \ge V_{\lsyse}^-(x,t)$. $V_{\lsys}$ must be bounded by $V_{\lsyse}^-$ and $V_{\lsyse}^+$, such that
    $
        V_{\lsys}(x,t) \in [V_{\lsyse}^-(x,t), V_{\lsyse}^+(x,t)].
    $
    \end{Lemma}
    \begin{proof}
        By definition $0 \in \eset$, thus $\stre_0 \triangleq 0 \in \stratsete$ which trivially implies $\tj_{\lsyse}(\tau; \stre_0) = \lintj(\tau)$. Therefore, taking the supremum or infimum over $\stratsete$ gives the desired bounds.
    \end{proof}

    In summary, it must be that for any $x \in \fset_c(t)$ and $t$,
    \begin{eqnarray}
        \vert V(x, t) - V_\lsys(x,t) \vert \le V_{\lsyse}^+(x, t) - V_{\lsyse}^-(x,t).
    \end{eqnarray}
    Note, by \cite{Evans84}, $V_{\lsyse}^+$ and $V_{\lsyse}^-$ are the viscosity solutions of HJ-PDE's with $H_{\lsyse}^+$ and $H_{\lsyse}^-$ given in \eqref{errhams}. Finally, since these HJ-PDE's have linear dynamics and convex $\tg$, when $H_{\lsyse}^\pm$ are convex, \cite{rublev2000generalized} certifies that
    \begin{eqnarray}
         V_{\lsyse}^+(x,t) =\hopf[\tgfn, H_{\lsys+\err}^+](x,t), \quad \text{ and } \quad V_{\lsyse}^-(x,t) =\hopf[\tgfn, H_{\lsys+\err}^-](x,t),
    \end{eqnarray}
    giving the desired relationship. Note, this bound may be guaranteed for concave $H_{\lsyse}^-$ and $H_{\lsyse}^+$ by reinitialization with sufficiently small $t$ (\cite{wang2019global}).

\end{proof}

\subsection{Proof of Corollary~\ref{cor:linbdtay}}

\begin{proof}
     Let $\ell$ be defined as in \eqref{taylorlin}, then by Taylor's theorem, $\delta^* = \max_{\mathcal{S} \times \ctrlset \times \distset \times \tint} \Vert [f - \ell ](x,u,d,t)\Vert \sim O(\Vert m - m_0 \Vert^2)$.  By Thm.~\ref{thm:conlinbd}, it must be that $V_\ell(x,t) \in [0, \epsilon^*]$ for any $x$ in a local region $\mathcal{S}$ (where $\epsilon^*$ is a function of the maximum error $\delta^*$ on $\mathcal{S}$), and thus there must exist an $\epsilon \in \mathbb{R}$ s.t. $V - V_\ell = \epsilon$. Considering this, one may note for any $\tau \in \tint$
     \begin{eqnarray*}
        \begin{aligned}
            \Vert \tj(\tau) - \lintj(\tau) \Vert 
            &= \left\Vert x + \int_t^s \sys(\tj(s), \csig(s), \strd[\csig](s))\:ds - \left(x + \int_t^s \lsys ( \lintj(s), \csig(s), \strd[\csig](s))\:ds \right) \right \Vert \\
            &= \left\Vert \int_t^s \sys(\tj(s), \csig(s), \strd[\csig](s)) - \sys ( \lintj(s), \csig(s), \strd[\csig](s)) + O(\Vert m - m_0 \Vert^2) \:ds \right \Vert \\
            &\le L_\lsys \int_t^s \Vert P \augtj(\tau) - \lintj(\tau) \Vert \: ds
        \end{aligned}
    \end{eqnarray*}
    and thus,
    \begin{eqnarray*}
        \begin{aligned}
            \lim_{m \to m_0} \Vert \tj(\tau) - \lintj(\tau) \Vert 
            &= \left\Vert \int_t^s \sys(\tj(s), \csig(s), \strd[\csig](s)) - \sys ( \lintj(s), \csig(s), \strd[\csig](s)) \:ds \right \Vert \\
            &\le L_\sys \int_t^s \Vert \tj(\tau) - \lintj(\tau) \Vert \: ds.
        \end{aligned}
    \end{eqnarray*} 
    This relation implies (\cite{boyce2021elementary}) that $ \lim_{m \to m_0} \Vert \tj(\tau) - \lintj(\tau) \Vert=0 \implies \lim_{m \to m_0} \lintj(\tau) = \tj(\tau)$. It then follows, 
    $$\lim_{m \to m_0} V_\ell(x,t) = \lim_{m \to m_0} \supod \infou J(\lintj(\tf)) = \supod \infou J(\tj(\tf)) = V(x,t).$$
\end{proof}

\subsection{Proof of Nonlinear Spectrum Augmentation Theorem \ref{thm:Vlam}}

\begin{proof} The proof of the $\augp$ boundary conditions follows from the fact that for initial condition $(x, \augp, t)$ at any time $\tau \in [t,\tf]$,
    % \begin{eqnarray*}
    % P \augtj(\tau) = \begin{cases}
    % \tj_\ell(\tau) & \text { if } \: \lambda = 0, \\
    % \tj(\tau)  & \text { if } \:\lambda = 1. \\
    % \end{cases}
    % \end{eqnarray*}
    \begin{eqnarray}
    P \augtj(\tau; x, \augp, t, \csig, \strd[u](\cdot)) = \begin{cases}
    \lintj(\tau; x, t, \csig, \strd[u](\cdot)) & \text { if } \: \augp = 0, \\
    \tj(\tau; x, t, \csig, \strd[u](\cdot))  & \text { if } \:\augp = 1. \\
    \end{cases}
    \label{augproj}
    \end{eqnarray}
    The proof of this fact for each case is analogous so we give that of $\augp=0$. Consider
    \begin{eqnarray*}
        \begin{aligned}
            \Vert P \augtj(\tau) - \lintj(\tau) \Vert 
            &= \left\Vert x + \int_t^s P\augsys(P\augtj(s), \csig(s), \strd[\csig](s))\:ds - \left(x + \int_t^s \lsys ( \lintj(s), \csig(s), \strd[\csig](s))\:ds \right) \right \Vert \\
            &= \left\Vert \int_t^s \lsys(P\augtj(s), \csig(s), \strd[\csig](s)) - \lsys ( \lintj(s), \csig(s), \strd[\csig](s))\:ds \right \Vert \\
            &\le L_\lsys \int_t^s \Vert P \augtj(\tau) - \lintj(\tau) \Vert \: ds
        \end{aligned}
    \end{eqnarray*}
    where $L_\lsys$ is a Lipschitz constant. This relation implies that $\Vert P \augtj(\tau) - \lintj(\tau) \Vert=0$; see a standard trajectory uniqueness proof for further details (\cite{boyce2021elementary}). 
    Note, the invariance of $\augp$ sufficed to fix the projection of $P \augsys$.
    With (\ref{augproj}), the value defined in (\ref{def:auggame}) is then trivially,
    \begin{eqnarray}
    V_\augp(x, \lambda, t) = \begin{cases} \supod \infou \minot \tgfn(\lintj(\tf))  & \text { if } \: \augp = 0, \\
    \supod \infou \minot \tgfn(\tj(\tf))  & \text { if } \:\augp = 1, \\
    \end{cases}
    \end{eqnarray}
    giving the boundary properites in (\ref{auggamebc}). Given the assumptions on the original game and that $\augsys$ is Lipschitz continuous w.r.t. $\augp$, it must be that $V_\augp$ is the solution of the HJ-PDE given in (\ref{HJPDE-Vlam}) and Lipschitz continuous w.r.t. $\augp$ (Thm. 2., \cite{evans1984differential}).
    % \square
\end{proof}

\subsection{Proof of Benchmark Property \ref{rem:benchprop}}

\begin{proof} Proof.
    Note the unidirectional influence ensures trajectories are ``decomposable'' \cite{Chen2015ExactDecomposition} s.t.
    \begin{eqnarray*}
        P_i \tj_{x,t}^{\csig, \dsig} (\tau) = \tj_{i, P_i x,t}^{\csig_i, \dsig_i} (\tau), \quad \tau \in [t,\tf]
    \end{eqnarray*}
    where $\tj_i$ is a trajectory of (\ref{pubsub2d}) and $\tj$ is a trajectory of (\ref{pubsubNd}). It follows that
    \begin{eqnarray*}
    \begin{aligned}
    V(x, t) = \sup_{\frak{d}} \inf_{\csig} \tgfn(\tj_{x,t}^{\csig, \frak{d}} (\tf))
    = \sup_{\{\frak{d}_i\}} \inf_{\{\csig_i\}} \sum^{N-1} J_{\tg,i}(\tj_{i, P_i x,t}^{\csig_i, \frak{d}_i} (\tf))
     = \sum^{N-1} V_i(P_i x, t).
    \end{aligned}
    \end{eqnarray*}
    It immediately follows that for $\tilde x \in \tilde{\mathcal{X}}$ and $i>0$,
    \begin{eqnarray*}
        V(\tilde x, t) = \sum^{N-1} V_i(P_i \tilde x, t) = (N-1) V_i(P_i \tilde x, t) \implies \bigg( V(\tilde x, t) = 0 \iff  V_i(P_i \tilde x, t) = 0 \bigg).
    \end{eqnarray*}
\end{proof}

\subsection{Training Details}\label{appendix:quad_dyn}

In this work, we introduced two new training programs as augmentations of the well-known existing approach proposed in \cite{bansal2020deepreach}, where the DeepReach software was proposed. For this reason, we adopted the parameters found in that work and our fork of the existing DeepReach software may be found \href{https://github.com/willsharpless/deepreach}{here}. As such, all neural nets used in this work consist of 3 layers with 512 neurons and sinusoidal activations. 

For the benchmark problem, we use the same training parameters for all variations of the system. 
For the baseline and the $V_\augp$ programs, we found that experiments with 300k iterations, a batch size of 60k and a learning rate of 5e-6 performed best on the 50-D case with our compute budget. For the decay program, we found that the experiments with 10k iterations, a batch size of 10k and a learning rate 1e-6 achieved the given results. For both LSS programs, the linear supervisors generated with or without Hopf data ultimately yielded equivalent accuracy but required five and ten minutes respectively; note, this run time is added to the final run time in all results.

In the quadrotor problem, we used a learning rate of 1e-5 for all variations of the system. 
For the baseline and the $V_\augp$ programs, we found that experiments with 100k iterations, and a batch size of 65k performed best. For both the linear and adaptive decay program, we used a batch size of 10k and 60k iterations for training. We further conducted a coarse parameter search for the LSS Decay program, and found that the experiments with $\lambda_K=0.6$, $\rho=0.1$, $\rho_g=0.2$ achieved the given results. For the LSS Decay Adp. method, we set $\mathcal{I}_{start}=10$ and $\mathcal{I}_{end} =1$.

\subsection{10-D Quadrotor Dynamics}\label{appendix:quad_dyn}
The quadrotor dynamics is given by
\begin{equation}
\begin{aligned}
    \dot{p_x}&= v_x \\
    \dot{p_y}&= v_y \\
    \dot{p_z}&= v_z \\
    \dot{\phi}&= - d_1 \phi + \omega_x \\
    \dot{\theta}&= - d_1 \theta + \omega_y \\
    \dot{v_x}&= g \ tan(\theta) \\
    \dot{v_y}&= g \ tan(\phi) \\
    \dot{v_z}&= u_3 \\
    \dot{\omega_x}&= - d_0 \phi + n_0 u_1 \\
    \dot{\omega_y}&= - d_0 \theta + n_0 u_2, \\
\end{aligned}
\end{equation}
where $(u_1, u_2, u_3) \in \left[-\frac{\pi}{4}, \frac{\pi}{4} \right]^2 \times \left[-1, 1 \right]$ is the control input and $(d_0, d_1, n_0)= (7,4,12)$ are constants. We consider states $\mathbb{X} \triangleq\left[-4, 4 \right]^2 \times \left[-2, 2 \right] \times \left[-1.5, 1.5 \right]^2 \times \left[-3, 3 \right]^2 \times \left[-2, 2 \right] \times \left[-6, 6 \right]^2$. 

\subsection{Adaptive Weighting Scheme}\label{appendix:adaptive_weight}
The overall idea of the adaptive weighting scheme is to gradually decrease the contribution of the supervision loss terms during training. To achieve this, we propose the following algorithm:

\begin{algorithm}[H]
\caption{Adaptive Weighting training}\label{alg:adap}
\begin{algorithmic}
\Require $N_k$, $\mathcal{I}_{end}$, $\mathcal{I}_{start}$
\Ensure $\mathcal{I}_{end} < \mathcal{I}_{start}$
\State 
\For{$k = 0,1,\ldots, K$}{
\State Compute importance of supervision loss terms $\mathcal{I}_k = \mathcal{I}_{end} \cdot exp\left\{ln(\frac{\mathcal{I}_{start}}{\mathcal{I}_{end}})*(1-\frac{k}{K})\right\}$
\State Compute loss terms $\losstermls(\theta)$ and $\losspde(\theta)$
\State Compute the relative weight $\lambda_k = 0.9\lambda_k + 0.1\mathcal{I}_k \frac{||\nabla_\lambda \losstermls(\theta)||}{||\nabla_\lambda \losspde(\theta)||}$
\State Compute overall loss $\loss_{LSS-A}(\theta; k) \triangleq \augp_k \mathbb{E}_{\ssp, \tint}[\losstermls(\theta)] + \mathbb{E}_{\ssp, \tint} [\losstermpde(\theta)]$
\State Step the optimizer}
\end{algorithmic}
\end{algorithm}
Empirically, choosing $\mathcal{I}_{start} = 10.0$ and $\mathcal{I}_{end} =1.0$ works well. 

% \begin{algorithm}
% 	\caption{Adaptive Weighting training} 
% 	\begin{algorithmic}[1]
% 		\For {$k = 1,2,\ldots, N_k$}
%                 \State Compute importance of supervision loss terms $\mathcal{I}_k = w_{end} \cdot exp\left\{ln(\frac{w_{start}}{w_{decay}})*(1-\frac{k}{N_k})\right\}$
%                 \State Compute loss terms $\losstermls(\theta)$ and $\losspde(\theta)$
% 				\State Compute the relative weight $\lambda_k = 0.9\lambda_k + 0.1\mathcal{I}_k \frac{||\nabla_\lambda \losstermls(\theta)||}{||\nabla_\lambda \losstermls(\theta)||}$
%                 \State Compute overall loss $\loss_{LSS-A}(\theta; k) \triangleq \augp_k \mathbb{E}_{\ssp, \tint}[\losstermls(\theta)] + \mathbb{E}_{\ssp, \tint} [\losstermpde(\theta)]$
%                 \State Step the optimizer
% 		\EndFor
% 	\end{algorithmic} 
% \end{algorithm}

\subsection{Metrics}\label{appendix:metrics}
\subsubsection{Probabilistic expansion}
We leverage method proposed in \cite{lin2023verification} to generate a probabilistic expansion from the learned value function using conformal prediction. To do this, we determine with high-confidence $\delta$, the maximum value for the empirically derived unsafe set under the learned policy, and then solve the corrected unsafe set $\hjrset_\delta = \left\{x: x \in \mathcal{X}, V_\theta(x,t)<\delta \right\}$. Then with $1-10^{-16}$ confidence, we can assert that 99.9\% of the states within the compliment $\hjrset_\delta^C$ are safe under the learned policy. 
% As one of the metrics, the volume of $\mathcal{S}$ is computed using 
% \begin{equation}
%     v_{\mathcal{S}}=\sum_{i=1}^{N} \frac{\mathbbm{1}(x_i \in \mathcal{S})}{N}, \ x_i \in Uniform(\mathcal{X}).
% \end{equation}
Finally, let the sampled volume of the recovered safe set be 
\begin{equation}
    \%v_{\mathcal{S}}=\mathbb{P}_{x \sim Uniform(\mathbb{X})}(x \in \hjrset_\delta^C).
\end{equation}
Note, the higher the volume the better the quality of learned solution.
\subsubsection{False positive rate \& False negative rate}
The false positive rate represents the likelihood of a random state $x\in \mathbb{X}$ being predicted as safe when it is unsafe under the learned policy, given by: 
\begin{equation}
    FP\%= \mathbb{P}_{x \sim Uniform(\mathbb{X})}(x \in \mathcal{S} \  \& \ \tgfn (\tj(t_f))<0).
\end{equation}
Similarly, the false negative rate is defined as the likelihood of a random state being classified as unsafe when it is indeed safe:
\begin{equation}
    FN\%= \mathbb{P}_{x \sim Uniform(\mathbb{X})}(x \in \mathcal{S}^C \  \& \ \tgfn (\tj(t_f))\geq0).
\end{equation}
Intuitively, $FP\%$ indicates the proportion of overly optimistic classifications while $FN\%$ represents the proportion of overly conservative classifications.

\end{document}

%% file: l4dc_2024/packages.tex
\usepackage[utf8]{inputenc} % allow utf-8 input
\usepackage[T1]{fontenc}    % use 8-bit T1 fonts
\usepackage{hyperref}       % hyperlinks
\usepackage{url}            % simple URL typesetting
\usepackage{booktabs}       % professional-quality tables
\usepackage{amsfonts}       % blackboard math symbols
\usepackage{nicefrac}       % compact symbols for 1/2, etc.
\usepackage{microtype}      % microtypography
\usepackage{xcolor}         % colors
\usepackage{wrapfig}
\usepackage{graphicx} 
\usepackage{graphicx} 
\usepackage{algorithm} 
\usepackage{algpseudocode} 
\usepackage{amsmath,amssymb}
\usepackage[font=small, labelfont=bf]{caption}

\usepackage{bbm}
\usepackage{soul}
\usepackage{microtype,wrapfig}
\usepackage{multicol}

\usepackage[export]{adjustbox}
\usepackage{makecell}
\usepackage{wrapfig}
\usepackage{bm}

\usepackage{lipsum}